\documentclass[leqno]{amsart}
\input
\input xy
\xyoption{all}
\usepackage{amssymb}
\newcommand{\rarrow}{\longrightarrow}
\def\U{{\mathfrak{A}}}
\def\A{\mathfrak{A}}
\def\G{\mathfrak{G}}

\def\LA{\bf LA}

\def\grL{{\bf LA}_{\Z_{\rm 2}}}

\def\LTS{\bf {LTS}}
\def\LTSfd{\bf {LTS}^{fd}}
\def\LAfd{\bf{LA}^{fd}_{\Z_{\rm 2}}}
\def\SSp{\bf{SymSp}^{(p)}}
\def\cscs{\bf {SymSp}_{csc}^{(p)}}
\def\cscg{\bf {LG}_{csc}^{\sigma}}

\def\LGi{\bf {LG}^{\sigma}}

\def\T{\mathcal T}
\def\R{\mathcal R}
\def\Lie{\mathfrak{d}}
\def\SSf{\mathfrak{S}}

\def\End{\mbox{End}}
\def\Hom{\mbox{Hom}}
\def\GrHom{\mbox{GrHom}}
\def\Kr{\mbox{\rm Ker}}
\def\Im{\mbox{\rm Im}}
\def\der{{\mathfrak{Der}}}
\def\inder{{\mathfrak{Inder}}}
\def\id{\mbox{\rm id}}
\def\span{{\rm span}}
\def\st{{\mathfrak Ste}}
\def\Cen{\mbox{\rm Center}}
\def\tan{{\rm T}}

\newcommand{\beq}{\begin{equation}}
\newcommand{\eeq}{\end{equation}}

\newtheorem*{theorema}{Theorem A}
\newtheorem*{theoremb}{Theorem B}
\newtheorem{lm}{Lemma}[section]

\newtheorem{thm}[lm]{Theorem}

\newtheorem{rem}[lm]{Remark}

\newtheorem{cor}[lm]{Corollary}

\newfam\msbfam
\font\tenmsb=msbm10 \textfont\msbfam=\tenmsb \font\sevenmsb=msbm7
\scriptfont\msbfam=\sevenmsb \font\fivemsb=msbm7
\scriptscriptfont\msbfam=\fivemsb

\def\to{\rightarrow}

\def\andd{\quad\mbox{and}\quad}

\def\Z{{\mathbb Z}} 

\begin{document}
\title{Imbedding of Lie triple systems into Lie algebras}
\author{O.\,N.\,Smirnov}
\date{January 11, 2005}

\begin{abstract}
We show that the category of Lie triple systems is equivalent to the category of $\Z_2$-graded Lie algebras $L=L_0\oplus L_1$ such that $L_1$ generates $L$ and the second graded cohomology group of $L$ with coefficients in any trivial module is zero. As a corollary we obtain an analogous result for symmetric spaces and Lie groups.
\end{abstract}

\maketitle

\section{Preliminaries.}
In this paper we will explicate connections between the category of Lie triple systems  and the category of Lie algebras.

Lie triple systems appeared initially in E.\,Cartan's study of symmetric spaces and totaly geodesic spaces~\cite{C}. The role played by Lie triple systems in the theory of symmetric spaces is parallel to that of Lie algebras in the theory of Lie groups: the tangent space at every point of a symmetric space has the structure of a Lie triple system. Moreover, the category of real finite-dimensional Lie triple systems is equivalent to the category of connected simply connected symmetric spaces with base point~\cite[Theorem~II.4.12]{L} or, alternatively, to the category of germs of symmetric spaces~\cite[Theorem I.2.9]{B}.

Other applications of Lie triple systems include representations of Jordan algebras~\cite{J1} and meson algebras~\cite{J2}, structurable  algebras and triple systems~\cite{F}, and Harish-Chandra modules over algebraic groups in positive characteristic~\cite{Ho2}.

One of the main tools in the structure theory of Lie triple systems is the universal imbedding of a Lie triple system into a Lie algebra constructed by N.\,Jacobson in~\cite{J1}. Its universal properties guarantee a functor $\A$ from the category of Lie triple systems to the category of Lie algebras. It is more natural to consider $\A$ as a functor to the category of $\Z_2$-graded Lie algebras. As such $\A$ is full, faithful, and constitutes a left adjoint to a natural ``forgetful'' functor. Our main result is the following description of the image of $\A$.
\begin{theorema} The category of Lie triple systems over a field is equivalent to the category of $\Z_2$-graded Lie algebras $L=L_0\oplus L_1$ such that
\begin{itemize}
\item[(i)] $L$ is generated by $L_1$ and \\[-0.3cm]
\item[(ii)] the second graded cohomology group $H^2_{\rm gr}(L,M)=0$ for every trivial module $M$ with the trivial grading $M=M_0$.
\end{itemize}
\end{theorema}

This paper is organized as follows. We recall necessary definitions and constructions in Section~\ref{sn2}.

In Section~\ref{sn3} we introduce a new construction of a universal enveloping Lie algebra $\U(T)$ for a Lie triple system $T$. It is done in the spirit of~\cite{AG} and~\cite{BS}. Our construction is more explicit than Jacobson's original construction and more convenient to work with.

Section~\ref{sn4} is concerned with the functor $\A$. To describe the image of $\A$ we develop a theory parallel to that of central extensions of perfect Lie algebras. We show that certain central extensions of $\U(T)$ are always trivial. To complete the proof of Theorem\,A we give a homological interpretation of this property.

In Section~\ref{sn5} we prove the following analog of Theorem A for Lie groups and symmetric spaces.

\begin{theoremb}
The category of connected simply-connected symmetric spaces with base point is equivalent to the category of Lie groups with involution $(G,\sigma)$ such that
\begin{itemize}
\item[(i)] $G$ is connected and simply-connected, \\[-0.3cm]
\item[(ii)] $G$ is generated by the space of symmetric elements $G_{\sigma}$, and \\[-0.3cm]
\item[(iii)] the invariant second cohomology group $H^2_{\rm inv}(G,\sigma)=0$.
\end{itemize}
\end{theoremb}

In our considerations the base ring is denoted by $k$ and all maps are assumed in to be $k$-linear.
Through the most of the paper $k$ is an arbitrary commutative ring. However, in our main result we need to assume that $k$ is a field to ensure the canonical connection between the second cohomology groups and central extensions.

\section{Lie triple systems and Lie algebras} \label{sn2}
In this section we recall necessary definitions, constructions, and results about Lie triple systems and the related Lie algebras. Details can be found in~\cite{J2}, \cite{H}, and~\cite{Ho}.
\subsection{Lie triple systems} \label{lts} A {\it Lie triple system} over a commutative ring $k$ is a $k$-module $T$ with a $k$-trilinear map $[\ ,\ ,\ ]:T\times T\times T\to T$ satisfying the identities
\begin{eqnarray}
&\left[a,a,b\right]=0 \label{lts1} \\
&\left[a,b,c\right]+[b,c,a]+[c,a,b]=0 \label{lts2}\\
&\left[a,b,[c,d,e]\right]=[[a,b,c],d,e]+[c,[a,b,d],e]+[c,d,[a,b,e]] \label{lts3}
\end{eqnarray}
for any $a,b,c,d,e\in T$. Note that the linearization of~(\ref{lts1}) wields
\beq
[a,b,c]+[b,c,a]=0. \label{lts4}
\eeq

A $k$-linear map $\alpha:T\to S$ between two Lie triple systems $T$ and $S$ is called a {\it homomorphism} of Lie triple systems if  $\alpha([a,b,c])=[\alpha(a),\alpha(b),\alpha(c)]$ for any $a,b,c\in T$.
The category of all Lie triple systems with the morphisms as above is denoted by $\LTS$.

Any Lie algebra $L$ considered with the trilinear product $[a,b,c]=[[a,b],c]$ is a Lie triple system.
More generally, any $k$-submodule $T$ of a Lie algebra $L$, which closed under $[[\ ,\ ],\ ]$, is a Lie triple system. Any Lie triple system can be realized in this way (see Sect.~\ref{sea}).

Since any Lie algebra homomorphism preserves the product $[a,b,c]$, one has a forgetful functor $\T:\LA\rarrow \LTS$ from the category $\LA$ of Lie algebras to $\LTS$: $\T L$ is a Lie algebra $L$ considered as a Lie triple system and $\T \varphi:\T L\to \T K$ is a Lie algebra homomorphism $\varphi:L\to K$ considered as a homomorphism of triple systems.

\subsection{Universal imbeddings of a Lie triple system}
If $T$ is a Lie triple system and $L$ is a Lie algebra, then a homomorphism $\varepsilon:T\to \T L$ is called an {\it imbedding} of the triple system $T$ into the algebra $L$. An imbedding $\nu:T\to \T U$ into a Lie algebra $U$ is called {\it universal} if for every imbedding $\varepsilon:T\to \T L$ there is a unique Lie algebra homomorphism $\varphi:U\to L$ such that the diagram below is commutative.
$$
\xymatrix{
T\ar[r]^{\nu}\ar[dr]_{\varepsilon}&\T U\ar[d]^{\T \varphi}\\
&\T L
}
$$
It was shown by Jacobson in~\cite{J1} that for every Lie triple system there exists a universal imbedding. Hence, the functor $\T:\LA\rarrow \LTS$ has a left adjoint.

Let us recall Jacobson's construction of a universal imbedding. For a Lie triple system $T$ let $L\langle T\rangle$ be a free Lie algebra on a $k$-module $T$ with the inclusion map $\iota:T\to L\langle T\rangle$ and let $I\langle T\rangle$ be the ideal of $L\langle T\rangle$ generated by the elements of the form $[[\iota(a),\iota(b)],\iota(c)]-\iota([a,b,c])$ for $a,b,c\in T$. Then the algebra $U=L\langle T\rangle/I\langle T\rangle$ together with the map $\nu:a\mapsto \iota(a)+I\langle T\rangle$ is a universal imbedding of $T$.

The fact that the map $\nu:T\to \T U$ is injective and other properties of $U$ follow from the next construction.

\subsection{The standard imbedding of a Lie triple system} \label{sea} An endomorphism $D\in \End_k(T)$ of a Lie triple system $T$ is said to be a {\it derivation} of $T$ if for any $a,b,c\in T$
\beq
D[a,b,c]=[Da,b,c]+[a,Db,c]+[a,b,Dc]. \label{der}
\eeq
It is easy to check that the set $\der(T)$ of all derivations of $T$ is a Lie subalgebra of $\End_k(T)$.

Setting $D_{a,b}c=[a,b,c]$ we obtain an endomorphism $D_{a,b}\in \End_k(T)$. Identity~(\ref{lts3}) implies that $D_{a,b}$ is a derivation of $T$.
Also,  equation~(\ref{der}) can be written in the form
\beq
[D,D_{a,b}]=D_{Da,b}+D_{a,Db} \label{inder}
\eeq
which implies that the set $\inder(T)=\span\{D_{a,b}:a,b\in T\}$ is an ideal of the Lie algebra $\der(T)$. The elements of $\inder(T)$ are called {\it inner derivations} of $T$.

The set $\st(T)=\inder(T)\oplus T$ with the product
\beq
[X+a,Y+b]=([X,Y]+D_{a,b})+(Xb-Ya), \label{st}
\eeq
where $X,Y\in\inder(T)$ and $a,b\in T$, is a Lie algebra. This algebra together with the canonical inclusion $T\to\inder(T)\oplus T$ is called the {\it standard imbedding} of the Lie triple system $T$~\cite{J1}.

\subsection{$\Z_2$-gradings} Recall that a Lie algebra $L$ is said to be {\it $\Z_2$-graded} if $L$ is a direct sum of a pair of $k$-submodules $L_0$ and $L_1$ such that $[L_i,L_j]\subseteq L_{i+j}$ for any $i,j\in\Z_2=\{0,1\}$. The  decomposition $L=L_0\oplus L_1$ is called a {\it $\Z_2$-grading} of $L$. If $L=L_0\oplus L_1$ and $K=K_0\oplus K_1$ are $\Z_2$-graded, a homomorphism $\varphi:L\to K$ is said to be {\it graded} provided that $\varphi(L_i)\subseteq K_i$ for any $i\in \Z_2$.

It follows immediately from~(\ref{st}) that the algebra $\st(T)$ is $\Z_2$-graded: $\st(T)_0=\inder(T)$ and $\st(T)_1=T$. Any universal imbedding $(U,\nu)$ of $T$ also possesses a $\Z_2$-grading. Indeed, the universal properties of $(U,\nu)$ imply that $\nu(T)$ generates $U$ and hence $U=[\nu(T),\nu(T)]+\nu(T)$. Moreover, since the canonical inclusion map $\iota:T\to\st(T)$ is injective and $[\iota(T),\iota(T)]\cap \iota(T)=0$, we conclude that $\nu:T\to U$ is injective and that $[\nu(T),\nu(T)]\cap \nu(T)=0$. Thus, the decomposition $U=[\nu(T),\nu(T)]\oplus\nu(T)$ is a $\Z_2$-grading of $U$.

Usually, the algebras $\st(T)$ and $U$ are considered as algebras with period-two automorphisms rather than as $\Z_2$-graded algebras. We prefer to consider $\Z_2$-gradings because it allows us to maintain a characteristic-free approach. Note also that these concepts are equivalent if the base ring $k$ has no 2-torsion.

\section{Alternative construction of universal imbedding} \label{sn3}

The aim of this section is to introduce a more explicit construction of the universal imbedding. Similar constructions were considered for the Lie algebras graded by root systems (see~\cite{ABG1}, \cite{ABG2}, and~\cite{BS}) and for certain $\Z$-graded algebras (see~\cite{AG}).

\subsection{Construction of $\U(T)$} Through the rest of this section we let $T$ be a Lie triple system over $k$. Our first step is to describe a central extension of the algebra $\inder(T)$.

Since $T$ is a module over the Lie algebra $\der(T)$, the exterior product $T\wedge T$ over $k$ is a $\der(T)$-module under the action
\beq
D\cdot(a\wedge b)=Da\wedge b+a\wedge Db \label{tmod}
\eeq
for $D\in \der(T)$ and $a,b\in T$. Moreover, formula (\ref{tmod}) and identity~(\ref{inder}) imply that the map $\lambda:T\wedge T\to \der(T)$ defined by $\lambda(a\wedge b)=D_{a,b}$ is a module homomorphism from $T\wedge T$ to the adjoint module $\der(T)$ and that $\Im(\lambda)=\inder(T)$.

To construct the required central extension we use the following version of ~\cite[Lemma~5.6]{BS}. This version was suggested by J.\,R.\,Faulkner.

\begin{lm} \label{mod}
Let $M$ be a module over a Lie algebra $L$ and let $\lambda:M\to L$ be an $L$-module homomorphism, i.e., $\lambda$ is a $k$-linear map such that
\beq
\lambda(l\cdot m)=[l,\lambda(m)] \label{lmod}
\eeq
for $l\in L$ and $m\in M$. Then
\begin{itemize}
\item[{\rm (i)}] $A(M)=\span\{\lambda(m)\cdot m:m\in M\}$ is a submodule of $M$ and $$\Im(\lambda)\cdot \Kr(\lambda)\subseteq A(M)\subseteq \Kr(\lambda);$$\\[-0.2in]
\item[{\rm (ii)}] the map $\mu:M/A(M)\to L$ induced by $\lambda$ is an $L$-module homomorphism;\\
\item[{\rm (iii)}] the quotient module $Q=M/A(M)$ forms a Lie algebra relative to the product
\beq
[p,q]=\mu(p)\cdot q \label{modprod}
\eeq
for $p,q\in Q$; and
\item[{\rm (iv)}] $\Im(\lambda)$ is a subalgebra of $L$ and $\mu:Q\to \Im(\lambda)$ is its central extension.
\end{itemize}
\end{lm}

\begin{proof}
Linearization of the expression $\lambda(m)\cdot m$ implies that elements of the form $\lambda(m)\cdot n+\lambda(n)\cdot m$ belong to $A(M)$. Then it follows from~(\ref{lmod}) that
\begin{equation*}
\begin{split}
l\cdot(\lambda(m)\cdot m)=[l,\lambda(m)]\cdot m+\lambda(m)\cdot(l\cdot m)=\lambda(l\cdot m)\cdot m+\lambda(m)\cdot(l\cdot m)\in A(M)
\end{split}
\end{equation*}
for any $l\in L$ and $m\in M$. Hence, $A(M)$ is a submodule of $M$. By~ (\ref{lmod}) one has $\lambda(\lambda(m)\cdot m)=[\lambda(m),\lambda(m)]=0$ and therefore $\lambda(A(M))=0$. Besides, for any $\lambda(m)\in \Im(\lambda)$ and any $n\in\Kr(\lambda)$ we have $\lambda(m)\cdot n=\lambda(m)\cdot n+\lambda(n)\cdot m\in A(M)$. Thus, (i) is established. The assertion~(ii) is immediate.

To prove~(iii) we note first that
\beq
\mu(p)\cdot p=0 \andd \mu(p)\cdot q+\mu(q)\cdot p=0\label{amodprod}
\eeq
for $p,q\in Q=M/A(M)$ and hence the product defined by~(\ref{modprod}) is anticommutative. Next, it follows from~(ii) and~(\ref{amodprod}) that
\begin{equation*}
\begin{split}
[[p,q],r]+[[q,r],p]+[[r,p],q]=\mu(\mu(p)\cdot q)\cdot r+
\mu(\mu(q)\cdot r)\cdot p+\mu(\mu(r)\cdot p)\cdot q\\
=[\mu(p),\mu(q)]\cdot r-\mu(p)\cdot(\mu(q)\cdot r)+\mu(q)\cdot(\mu(p)\cdot r)=0
\end{split}
\end{equation*}
for $p,q,r\in Q$. Thus $Q$ is a Lie algebra.

Finally, to establish (iv) we note that the equality  $\mu([p,q])=\mu(\mu(p)\cdot q)=[\mu(p),\mu(q)]$ follows from (ii) and~(\ref{modprod}). Hence,  $\Im(\lambda)=\Im(\mu)$ is a subalgebra of $L$ and $\mu$ is a homomorphism. Formula~(\ref{modprod}) implies also that $\Kr(\mu)\subseteq \Cen(Q)$. Thus, $\mu:Q\to \Im(\mu)=\Im(\lambda)$ is a central extension.
\end{proof}

Let us return to the map $\lambda:T\wedge T\to \der(T)$, $\lambda(a\wedge b)=D_{a,b}$. It is easy to see that $A(T\wedge T)=\span\{\lambda(x)\cdot x:x\in T\wedge T\}$ is spanned by the elements of the form
\begin{eqnarray}
&D_{a,b}(a\wedge b) \label{i1}\\
&D_{a,b}(c\wedge d)+D_{c,d}(a\wedge b) \label{i2}
\end{eqnarray}
for $a,b,c,d\in T$. Note that if $k$ has no 2-torsion, then~(\ref{i1}) is a special case of~(\ref{i2}).

We let $\langle T,T\rangle$ denote the quotient module $(T\wedge T)/A(T\wedge T)$ and let $\langle a,b\rangle$ denote the coset containing $a\wedge b$. Lemma~\ref{mod} implies
\begin{cor} \label{TT}
The $k$-module $\langle T,T\rangle$ is a Lie algebra relative to the product
\beq
[\langle a,b\rangle,\langle c,d\rangle]=\langle D_{a,b}c,d\rangle+\langle c,D_{a,b}d\rangle .\label{TTprod}
\eeq
Moreover, the map $\mu:\langle T,T\rangle\to\inder(T)$, $\mu(\langle a,b\rangle)=D_{a,b}$, is a central extension.
\end{cor}

To lift $\mu:\langle T,T\rangle\to\inder(T)$ to the extension of $\st(T)=\inder(T)\oplus T$ the following axillary lemma is useful. Its verification is straightforward and therefore is omitted.

\begin{lm}
Let $M$ be a module over a Lie algebra $L$ and let
$\langle \ ,\ \rangle:M\wedge M\to L$ be an $L$-module homomorphism, i.e.,
\beq
[x,\langle m,n\rangle]=\langle x\cdot m,n\rangle+\langle m,x\cdot n\rangle
\label{wed1}
\eeq
for $x\in L$ and $m,n\in M$. Assume also that
\beq
\langle m,n\rangle\cdot k+\langle n,k\rangle\cdot m+
\langle k,m\rangle\cdot n=0 \label{wed2}
\eeq
for every $m,n,k\in M$. Then the space $L\oplus M$ forms a $\Z_2$-graded Lie algebra relative to the product
\begin{equation*}
[x+m,y+n]=([x,y]+\langle m,n\rangle)+(x\cdot n-y\cdot m)
\end{equation*}
for $x,y\in L$ and $m,n\in M$.
\end{lm}

\begin{cor} \label{TTT} The space $\U (T)=\langle T,T\rangle\oplus T$ with the product
\beq
[X+a,Y+b]=([X,Y]+\langle a,b\rangle)+(\mu(X)a-\mu(Y)b), \label{Lprod}
\eeq
where $X,Y\in \langle T,T\rangle$ and $a,b\in T$, is a $\Z_2$-graded Lie algebra. Moreover, the map $\upsilon:\U (T)\to \st(T)$, defined by $\upsilon(X+a)=\mu(X)+a$, is a graded central extension.
\end{cor}

\begin{proof} First we note that $T$ is a $\langle T,T\rangle$-module relative to the representation $\mu:\langle T,T\rangle\to \inder(T)\subseteq \End_k(T)$. Besides,~(\ref{wed1}) follows from~(\ref{TTprod}) and~(\ref{wed2}) follows from~(\ref{lts2}).

Finally, comparing~(\ref{st}) and~(\ref{Lprod}) one can see that $\upsilon$ is a homomorphism. Moreover, $\upsilon$ is a central extension because $\mu$ is.
\end{proof}

\subsection{Universal property of $\U(T)$}

Here we consider $T$ as a subspace of $\U (T)=\langle T,T\rangle\oplus T$. It follows from~(\ref{Lprod}) that $T$ generates $\U (T)$.

\begin{thm} \label{Tuni} For every Lie triple system $T$ the canonical inclusion map $\iota_T:T\to\U (T)=\langle T,T\rangle\oplus T$ is a universal imbedding.
\end{thm}

\begin{proof} Let $L$ be a Lie algebra and let $\varepsilon:T\to \T L$ be an imbedding, i.e., for every $a,b,c\in T$
\beq
\varepsilon([a,b,c])=[[\varepsilon(a),\varepsilon(b)],\varepsilon(c)]. \label{a1}
\eeq
We need to construct a Lie algebra homomorphism $\varphi:\U (T)\to L$ extending $\varepsilon$ and prove that such $\varphi$ is unique.

First, we claim that there exists a well-defined map $\langle \varepsilon,\varepsilon \rangle:\langle T,T\rangle\to L$ such that
\beq
\langle \varepsilon,\varepsilon \rangle(\langle a,b\rangle)=[\varepsilon(a),\varepsilon(b)] \label{a2}
\eeq
for every $a,b\in T.$ It suffices to show that the $k$-submodule $A(T\wedge T)$ of $\langle T,T\rangle$ spanned by the elements~(\ref{i1}) and~(\ref{i2}) is in the kernel of the map $\zeta:T\wedge T\to L$ defined by $\zeta(a\wedge b)=[\varepsilon(a),\varepsilon(b)]$. This is indeed the case, since for any $a,b,c,d\in T$ we have
\begin{equation*}
\begin{split}
\zeta([a,b,a]\wedge b+a\wedge [a,b,b])=
[\varepsilon([a,b,a]),\varepsilon(b)]+[\varepsilon(a),\varepsilon([a,b,b])]\\
=[[[\varepsilon(a),\varepsilon(b)],\varepsilon(a)],\varepsilon(b)]+
[\varepsilon(a),[[\varepsilon(a),\varepsilon(b)],\varepsilon(b)]]\\
=[[\varepsilon(a),\varepsilon(b)],[\varepsilon(a),\varepsilon(b)]]=0 \ \andd\\
\end{split}
\end{equation*}
\begin{equation*}
\begin{split}
\zeta([a,b,c]\wedge d+c\wedge [a,b,d]+[c,d,a]\wedge b+a\wedge [c,d,b])\\
=[[[\varepsilon(a),\varepsilon(b)],\varepsilon(c)],\varepsilon(d)]
+[\varepsilon(c),[[\varepsilon(a),\varepsilon(b)],\varepsilon(d)]]+
[[[\varepsilon(c),\varepsilon(d)],\varepsilon(a)],\varepsilon(b)]\\
+[\varepsilon(a),[[\varepsilon(c),\varepsilon(d)],\varepsilon(b)]]\\
=[[\varepsilon(a),\varepsilon(b)],[\varepsilon(c),\varepsilon(d)]]+
[[\varepsilon(c),\varepsilon(d)],[\varepsilon(a),\varepsilon(b)]]
=0.
\end{split}
\end{equation*}

Since $T$ generates the Lie algebra $\U (T)$, the identities~(\ref{a1}) and~(\ref{a2}) imply that the map $\varphi$, defined by
\beq \varphi(X+a)=\langle \varepsilon,\varepsilon \rangle(X)+\varepsilon(a)  \label{fi}
\eeq
for $X\in \langle T,T\rangle$ and $a\in T$, is a Lie algebra homomorphism. By construction, $\varphi|_T=\varepsilon$. Since $T$ generates $\U (T)$, there is only one such homomorphism $\varphi$.
\end{proof}

\section{A category of Lie algebras equivalent to $\LTS$} \label{sn4}

The goal of this section is to determine a category of Lie algebras equivalent to $\LTS$.

\subsection{The functor $\A:\LTS\rarrow \grL$} \label{GrFunctor}
The universal property of $\U(T)$ guarantees the existence of a left adjoint to the forgetful functor $\T:\LA\rarrow \LTS$ discussed in Section~\ref{lts}. The functor $\A:\LTS\rarrow \LA$ sends every system $T$ to the algebra $\U(T)$ and every morphism $\alpha:T\to S$ to the morphism $\varphi:\U(T)\to\U (S)$ defined by~(\ref{a2}) and~(\ref{fi}) for the imbedding $\varepsilon=\iota_S\circ\alpha$.
These formulas imply also that the functor $\A$ is faithful.

In fact, it is advantageous to consider $\A$ as a functor to the category $\grL$ of $\Z_2$-graded Lie algebras with graded homomorphisms, because as such $\A$ is both full and faithful. Since $\grL$ is a subcategory of $\LA$ we use the same notation for the functor $\A:\LTS\rarrow \grL$.

Furthermore, $\A:\LTS\rarrow \grL$ has a right adjoint $\R$ which can be described as follows: for a $\Z_2$-graded Lie algebra $L=L_0\oplus L_1$ the odd component $L_1$ is a Lie triple system relative to the ternary product $[a,b,c]=[[a,b],c]$, so we set $\R L=L_1$. If $\varphi:L\to K$ is a graded homomorphism, then the restriction of $\varphi$ on $L_1$ gives us a homomorphism of Lie triple systems which is taken to be $\R\varphi$. Theorem~\ref{Tuni} implies that $\iota_T:T\to \R\U(T)\subseteq \U(T)$ is a universal arrow from $T$ to $\R$. For convenience of references we restate Theorem~\ref{Tuni} for $\Z_2$-graded algebras.

\begin{thm} \label{Tuniv}
Let $L=L_{0}\oplus L_{1}$ be a $\Z_2$-graded Lie algebra and let $T$ be a Lie triple system. For every homomorphism $\alpha:T\to L_1$ there is a unique graded homomorphism $\widehat\alpha:\U (T)\to L$ extending $\alpha$. This homomorphism is defined by the formula: $\widehat\alpha(\sum_i \langle t_i,t'_i\rangle+t)=\sum_i [\alpha(t_i),\alpha(t'_i)]+\alpha(t)$.
\end{thm}

Since $\A:\LTS\rarrow \grL$ is full and faithful, the category $\LTS$ is equivalent to the full subcategory of $\grL$ whose objects are the Lie algebras isomorphic to the algebras $\U(T)$. To give an intrinsic description of this subcategory we need to characterize the Lie algebras of the form $\U (T)$.

\subsection{Central 0-extensions} Our first goal is to show that the central extension $\upsilon:\U (T)\to\st(T)$ is universal in a certain class of extensions. The model for our treatment is the theory of central extensions of perfect Lie algebras (e.g. see~\cite[Sect.\,7.9]{W}). Throughout this section all Lie algebras and homomorphisms are assumed to be $\Z_2$-graded.

A homomorphism $\varphi: K\to L$ is called a 0-{\it extension} of $L$ if $\varphi$ is surjective and $\Kr(\varphi)\subseteq K_{0}$. Thus an epimorphism $\varphi:K\to L$ is a central 0-extension if $\Kr(\varphi)\subseteq K_{0}\cap\Cen(K)$.
Naturally, a central 0-extension $\upsilon:U\to L$ is termed {\it universal} if for every central 0-extension $\varphi:K\to L$ there is a unique homomorphism $\psi:U\to K$ such that the diagram
$$
\xymatrix{
U\ar[r]^{\upsilon}\ar[d]_{\psi}&L\\
K\ar[ur]_{\varphi}&
}
$$
is commutative. It follows that a universal central 0-extension of any Lie algebra is defined uniquely up to an isomorphism.

\begin{thm} \label{TuLT} If a $\Z_2$-graded Lie algebra $L$ is generated by $L_1$, then $\upsilon:\U(L_1)\to L$ is a universal central 0-extension of $L$ for the map $\upsilon=\widehat{(\mbox{\rm id}_{L_1})}$ as in Theorem~\ref{Tuniv}.
\end{thm}

\begin{proof}
Assume that $\varphi: K\to L$ is a central 0-extension of $L$. Since $\varphi$ is an epimorphism and $\Kr(\varphi)\subseteq K_0$, the Lie triple system homomorphism $\R\varphi:K_1\to L_1$ is invertible. It follows then from Theorem~\ref{Tuniv} that for the map $(\R\varphi)^{-1}:L_1 \to K_1$ there exists a unique Lie algebra homomorphism $\psi:\U(L_1)\to K$ such that $\R \psi=(\R\varphi)^{-1}$ or $\R\varphi\circ\R\psi=\mbox{\rm id}_{L_1}$. Since $\R\upsilon=\mbox{\rm id}_{L_1}$, one has $\R\varphi\circ\R\psi=\R\upsilon$. This equality implies $\varphi\circ\psi=\upsilon$ since $\U(L_1)$ is generated by $L_1$.
\end{proof}

\begin{cor} \label{0dis} If a $\Z_2$-graded Lie algebra $L$ is generated by $L_1$, then $L\simeq \langle L_1,L_1\rangle/I\oplus L_1$ where $I$ is a central ideal of $L$ contained in $\langle L_1,L_1\rangle$.
\end{cor}

\begin{cor} \label{TuLTC} For every Lie triple system $T$ the algebra $\U(T)$ is a universal central 0-extension of $\st(T)$.
\end{cor}

\begin{proof}
Since $\st(T)_1= T$, the proof follows.
\end{proof}

\begin{rem} {\rm
Theorem~\ref{TuLT} implies that if $L$ is generated by $L_1$ then any universal central 0-extension $U=U_0\oplus U_1$ of $L$ is generated by $U_1$. Moreover, $L$ is generated by $L_1$ if and only if there exists a universal central 0-extension $U$ of $L$ such that $U$ is generated by $U_1$. Unlike the case with perfect Lie algebras, mere existence of a universal central 0-extension of $L$ does not necessarily imply that $L$ is generated by $L_1$. For example, a non-zero perfect Lie algebra $L$ with the trivial grading $L=L_0$ is not generated by $L_1$. However, a universal central extension of $L$ exists and it is also a universal central 0-extension of $L$. }
\end{rem}

Our next step is to prove a recognition theorem for universal central 0-extensions similar to~\cite[Recognition Criterion 7.9.4]{W} for perfect Lie algebras. Here we say that a $\Z_2$-graded Lie algebra $L=L_{0}\oplus L_{1}$ is {\it 0-centrally   closed}, if every central 0-extension of $L$ splits.

\begin{thm} \label{TrecU} Assume that $\upsilon:U\to L$ is a central 0-extension of a Lie algebra $L$ generated by $L_1$. This extension is universal if and only if $U$ is generated by $U_1$ and 0-centrally closed.
\end{thm}

\begin{proof}
Assume first that $\upsilon:U\to L$ is universal. We already mentioned that in this case $U$ is generated by $U_1$.

Let $\varphi:K\to U$ be a central 0-extension of $U$. Then $K'=[K_{1},K_{1}]+K_1$  is a graded subalgebra of $K$, the inclusion map $\iota':K'\to K$ and the composition $\upsilon\circ\varphi\circ\iota':K'\to L$ are a graded homomorphisms.
We claim that $\upsilon\circ\varphi\circ\iota':K'\to L$ is a central 0-extension of $L$.

The inclusions $\Kr(\upsilon)\subseteq U_{0}$ and $\Kr(\varphi)\subseteq K_{0} $ imply that $\Kr(\upsilon\circ\varphi)=\{k\in K:\varphi(k)\in \Kr(\upsilon)\}\subseteq \{k\in K:\varphi(k)\in U_{0} \}\subseteq K_{0} $.
This together with the fact that $\upsilon\circ\varphi$ is an epimorphism imply that $(\upsilon\circ\varphi\circ\iota')(K_{1})=L_{1}$.
Furthermore, $L_{0}=[L_{1},L_{1}]=
(\upsilon\circ\varphi\circ\iota')([K_{1},K_{1}])=
(\upsilon\circ\varphi\circ\iota')(K'_{0})$.
Thus, $\upsilon\circ\varphi\circ\iota'$ is an epimorphism.

In addition, $\Kr(\upsilon\circ\varphi\circ\iota')\subseteq K'\cap K_{0} \subseteq K'_{0} $. Moreover, $\varphi(k)\in \Kr(\upsilon)\subseteq \Cen(U)$ for every $k\in\Kr(\upsilon\circ\varphi\circ\iota')$ and therefore we have $[k,K_{1}]\in \Kr(\varphi)\cap K_{1} =0$. This implies that $\Kr(\upsilon\circ\varphi\circ\iota')\subseteq \Cen(K')$ justifying our claim that $\upsilon\circ\varphi\circ\iota'$ is a central 0-extension.

For this extension there exists a map $\psi:U\to K'$ such that $(\upsilon\circ\varphi\circ\iota')\circ\psi=\upsilon$. However, for the extension $\upsilon:U\to L$ the identity map $\mbox{\rm id}_U$ has the property $\upsilon\circ\mbox{\rm id}_U=\upsilon$. So $\varphi\circ\iota'\circ\psi=\mbox{\rm id}_U$ and therefore $U$ is 0-centrally closed.

Conversely, assume that $U$ is generated by $U_1$ and is 0-centrally  closed. Let $\varphi:K\to L$ be a central 0-extension of $L$.

To construct a map from $U$ to $K$ we will show that the pullback of $\varphi$ along $\upsilon$ is a central $0$-extension of $U$. The direct sum of algebras $K\oplus U$ is endowed with a $\Z_2$-grading: $(K\oplus U)_i=K_i\oplus U_i$ for $i\in\Z_2$. It is easy to see that the set $A=\{k+u\in K\oplus U :\varphi(k)=\upsilon(u)\}$ is a graded subalgebra of $K\oplus U$ and that the projections $\pi_K:A\to K$, $\pi_K(k+u)=k$, and $\pi_U:A\to U$, $\pi_U(k+u)=u$, are graded epimorphisms. The definition of $A$ implies that
\beq
\varphi\circ \pi_K=\upsilon\circ\pi_U. \label{fu}
\eeq

Note that $\pi_U:A\to U$ is a central 0-extension of $U$ because $\Kr(\pi_U)=\{k+0\in K\oplus U:k\in\Kr(\varphi)\}\subseteq \Cen(A)\cap A_{0} $. Thus there is a splitting morphism $\psi:U\to A$ with $\pi_U\circ \psi=\mbox{\rm id}_U$. It follows then from~(\ref{fu}) that $\varphi\circ (\pi_K\circ\psi)=\upsilon\circ\pi_U\circ\psi=\upsilon$.

To show that $\pi_K\circ\psi$ is unique with this property assume that $\varphi\circ\mu=0$ for a homomorphism $\mu:U\to K$. One has $\mu(U)\subseteq \Kr(\varphi)\subseteq K_{0} $ and hence $\mu(U_{1})=0$. It follows that $\mu=0$ because $U$ is generated by $U_1$.
\end{proof}
\begin{cor}
A $\Z_2$-graded Lie algebra $L$ is isomorphic to $\U(T)$ for some Lie triple system $T$ if and only if $L$ is generated by $L_1$ and is 0-centrally closed.
\end{cor}
\begin{proof} If $L\simeq \U(T)$, then $L$ is generated by $L_1$. It follows from Corollary~\ref{TuLTC} and Theorem~\ref{TrecU} that $\U(T)$ is 0-centrally closed.

Conversely, assume that $L$ is generated by $L_1$ and 0-centrally closed. Then the central 0-extension $\upsilon:\U(L_1)\to L$ described in Theorem~\ref{TuLT} splits. The splitting map $\theta:L\to\U(L_1)$ is injective since $\upsilon\circ\theta=\id_L$.

We claim that $\theta$ is also bijective. Indeed, since $\upsilon|_{L_1}=\id_{L_1}$ for $l_1,k_1\in L_1$ in $\U(L)$ one has $\langle l_1,k_1\rangle=[l_1,k_1]=[\theta(\upsilon(l_1)),\theta(\upsilon(k_1))]=
\theta([\upsilon(l_1),\upsilon(k_1)])$. Thus, $L\simeq \U(L_1)$.
\end{proof}
\begin{cor} \label{cat1}
The category $\LTS$ of Lie triple systems is equivalent to the category of $0$-centrally closed $\Z_2$-graded Lie algebras generated by their odd graded components.
\end{cor}
\subsection{Homological characterization} \label{homi} The condition on extensions in Corollary~\ref{cat1} has a natural homological characterization which we will find in this section.

If $k$ is a field of characteristic different from 2, Harris developed a cohomology theory for the Lie algebras with involutions and showed that the central extensions in this category of algebras are in a bijective correspondence with the elements of the second cohomology group which preserve the involution~\cite{H}.

For our purposes we need a graded version of the classical cohomology theory. Here we present an adaptation of the classical notions and arguments to the graded algebras. Specifically, we consider a graded version of the second cohomology group $H^2(L,M)$ for a graded module $M$ over a $\Z_2$-graded Lie algebra $L$ and its relation to the graded central extensions of $L$.

Recall that the cohomology groups $H^*(L,M)$ of a Lie algebra $L$ with coefficients in an $L$-module $M$ can be defined as the cohomology groups of the Chevalley-Eilenberg cochain complex
$$
...\longrightarrow\Hom _k(\wedge^n L,M)\stackrel{\delta}{\longrightarrow}
\Hom _k(\wedge^{n+1} L,M)\longrightarrow ...
$$
where the coboundary map $\delta$ is given by (see~\cite{W} for details).
\beq\label{cobo}
\begin{split}
\delta\! f(x_1,x_2,...,x_{n+1})&=\sum_{i=1}^{n+1}(-1)^{i+1}x_i\cdot f(x_1,...,\hat x_i,...,x_{n+1})\\
&+\sum_{1\leq i<j\leq n+1}(-1)^{i+j}f([x_i,x_j],x_1,...,\hat x_i,...,\hat x_j,...,x_{n+1}).
\end{split}
\eeq
If $L$ is a $\Z_2$-graded algebra and $M$ is a graded module, one can introduce a graded version of the complex above. Let $\GrHom_k(\wedge^n L,M)$ be the set of elements $f$ from $\Hom _k(\wedge^n L,M)$ such that $f(L_{i_1},L_{i_2},...,L_{i_n})\subseteq M_{i_1+i_2+...+i_n}$ for every $i_1,i_2,...,i_n\in \Z_2$. It follows from~(\ref{cobo}) that $\delta(\GrHom_k(\wedge^n L,M))\subseteq \GrHom_k(\wedge^{n+1} L,M)$. Hence one has the complex
$$
...\longrightarrow\GrHom _k(\wedge^n L,M)\stackrel{\delta}{\longrightarrow}
\GrHom _k(\wedge^{n+1} L,M)\longrightarrow ...\quad .
$$
The cohomology groups of this complex are defined to be the {\it cohomology groups $H^*_{\rm gr}(L,M)$} of the graded Lie algebra $L$ with coefficients in the graded module $M$. In fact, each $H^n_{\rm gr}(L,M)$ can be identified with a subgroup of $H^n(L,M)$.

\begin{thm} \label{homolal}
If for a $\Z_2$-graded algebra $L=L_0\oplus L_1$ every central 0-extension splits, then $H^2_{\rm gr}(L,M)=0$ for every trivial module $M$ such that $M_1=0$.

Moreover, if $L_0$ is a projective $k$-module, then the converse is true.
\end{thm}

\begin{proof}
Assume that every central 0-extension of $L$ splits, $M$ is a trivial $L$-module, and $M_1=0$.

It is well known that any 2-cocycle $\sigma\in Z^2(L,M)$ defines a central extension $\varphi:L\ltimes_{\sigma} M\to L$. The algebra $L\ltimes_{\sigma} M$ is the direct sum of $k$-modules $L\oplus M$ with the product
\begin{equation}
\left[x+m,y+n\right]=[x,y]+\sigma(x,y) \label{sig}
\end{equation}
for $x+m,y+n\in L\oplus M$; and $\varphi:x+m\mapsto x$ is the projection map.

It is easy to see that if $\sigma\in Z^2_{\rm gr}(L,M)$, then the decomposition $L\ltimes_{\sigma} M=(L_0\oplus M_0)\oplus L_1$ is a $\Z_2$-grading and therefore $\varphi$ is a 0-extension. By our assumption there exists a splitting homomorphism $\tau:L\to L\ltimes_{\sigma} M$, that is $\tau(x)=x+\tau_M(x)$ for some linear map $\tau_M:L\to M$. It follows from~(\ref{sig}) that for every $x,y\in L$
\begin{equation*}
\begin{split}
[x,y]+\tau_M([x,y])=\tau([x,y])=[\tau(x),\tau(y)]=[x+\tau_M(x),y+\tau_M(y)]\\
=[x,y]+\sigma(x,y).
\end{split}
\end{equation*}
Thus $\sigma(x,y)=\tau_M([x,y])$ for any $\sigma\in Z^2_{\rm gr}(L,M)$ and therefore $H^2_{\rm gr}(L,M)=0$.

Assume now that $L_0$ is a projective $k$-module and that $H^2_{\rm gr}(L,M)=0$ for every trivial module $M$ such that $M_1=0$. Let $\varphi:K\to L$ be a central 0-extension of $L$. Since $\Kr(\varphi)\subseteq K_0$ and $L_0$ is projective there is a linear map $\eta:L\to K$ such that $\eta$ preserves the grading and $\varphi\circ \eta=\id_L$.

We consider $\Kr(\varphi)$ as a trivial $L$-module with the grading $\Kr(\varphi)=\Kr(\varphi)_0$. It is known that the map $\sigma:L\times L\to \Kr(\varphi)$ defined by $\sigma(x,y)=[\eta(x),\eta(y)]-\eta([x,y])$ is a 2-cocycle. Moreover, for every $x_i\in L_i$ and $y_j\in L_j$ we have $\sigma(x_i,y_j)=[\eta(x_i),\eta(y_j)]-\eta([x_i,y_j])\subseteq M\cap K_{i+j}\subseteq M_{i+j}$. Hence, $\sigma\in Z^2_{\rm gr}(L,\Kr(\varphi))$.

Since $H^2_{\rm gr}(L,\Kr(\varphi))=0$, there is a linear map $\tau:L\to \Kr(\varphi)$
such that $\sigma(x,y)=\tau([x,y])$. Then the map $\psi=\eta+\tau$ is a Lie algebra homomorphism $\psi:L\to K$. Indeed,
\begin{equation*}
\begin{split}
\psi([x,y]))=\eta([x,y])+\tau([x,y])=[\eta(x),\eta(y)]=
[\eta(x)+\tau(x),\eta(y)+\tau(y)]\\
=[\psi(x),\psi(y)].
\end{split}
\end{equation*}
Besides, $\varphi\circ\psi=\id_L$ since $\varphi\circ\tau=0$. Thus $\psi$ is a splitting map for $\varphi$.
\end{proof}

\begin{rem} {\rm
Theorem~\ref{homolal} and Corollary~\ref{cat1} furnish a proof for Theorem~A.}
\end{rem}

\section{Lie groups with involution and symmetric spaces} \label{sn5}

In conclusion we show that the imbedding of Lie triple systems into graded Lie algebras, described in Theorem A, can be lifted to an imbedding of symmetric spaces to Lie groups with involution.

To this end we consider the category $\LGi$ of Lie groups with involution. An  object of $\LGi$ is a pair $(G,\sigma_G)$ where $G$ is a real Lie group and $\sigma_G$ is an involutive automorphism of $G$, called an {\it involution} for short. We often write $\sigma$ for $\sigma_G$ , when no confusion can arise. A morphism from $(G,\sigma_G)$ to $(H,\sigma_H)$ is a Lie group homomorphism $\varphi:G\to H$ such that $\varphi\circ\sigma_G=\sigma_H\circ\varphi$.

Let $\Lie$ be the functor that sends a Lie group $G$ to its Lie algebra $g=\Lie G$. Then for every object $(G,\sigma)$ from $\LGi$ the involutive automorphism $\Lie \sigma$ of the Lie algebra $g=\Lie  G$ defines a $\Z_2$-grading $g=g_0\oplus g_1$ where $g_i=\{x\in g:x^{\Lie\sigma}=(-1)^ix\}$. Moreover, it is easy to see that $\Lie$ can be considered as a functor from $\LGi$ to the category of real $\Z_2$-graded Lie algebras $\grL$.

Next, we recall that a real manifold $M$ is termed
a {\it symmetric space} if for any $x\in M$ there is a differentiable map $S_x:M\to M$ such that
\begin{itemize}
\item[(i)]$S_x^2=\id_M$, \\
\item[(ii)] $S_xS_yS_x=S_{S_xy}$ for any $y\in M$, and \\
\item[(iii)] $x$ is an isolated fixed point of $S_x$.
\end{itemize}
The tangent space $\tan_p(M)$ at a point $p$ is a Lie triple system under the product $[X,Y,Z]=R(X,Y)Z$ where $R:M\otimes M\to M$ is the curvature tensor of $M$. In fact, the map $(M,p)\mapsto \tan_p(M)$ rises to a functor $\Lie_{\rm lts}$ from  the category $\SSp$ of pointed symmetric spaces to the category $\LTS$ of Lie triple systems.

On the other hand, for an object $(G,\sigma)$ from $\LGi$ the {\it space of symmetric elements} $G_{\sigma}=\{g\sigma(g)^{-1}:g\in G\}$ forms a symmetric space under the operation $S_x(y)=xy^{-1}x$. Furthermore, the restriction of any morphism $\varphi:(G,\sigma_G)\to (H,\sigma_H)$ gives a morphism $\varphi|_{G_{\sigma}}:(G_{\sigma},e_G)\to (H_{\sigma},e_H)$ of pointed symmetric symmetric spaces where $e_G$ and $e_H$ are the identity elements of $G$ and $H$ respectively. The functor thus obtained is denoted by $\SSf$.

These three functors together with the forgetful functor $\R:\grL\to \LTS$ introduced in ~\ref{GrFunctor} form the following diagram
$$
\xymatrix{
\LGi\ar[r]^{\SSf\quad}\ar[d]_{\Lie}&\SSp\ar[d]^{\Lie_{\rm lts}}\\
\grL\ar[r]^{\R}&\LTS
}
$$ which is commutative by Proposition IV.4.4 from~\cite{L}.

To invert the vertical arrows in this diagram one has to restrict them to appropriate subcategories. Namely, $\Lie$ provides an equivalence between the category $\cscg$ of connected simply connected Lie groups with involution and the category $\LAfd$ of finite-dimensional graded Lie algebras; while $\Lie_{\rm lts}$ furnishes the equivalence between the category $\cscs$ of connected simply connected symmetric spaces with base point and the category $\LTSfd$ of finite-dimensional Lie triple systems~\cite[Theorem~II.4.12]{L}.

It follows from Theorem A that the functor $\A$ is an imbedding of $\LTSfd$ into $\LAfd$, so composing $\A$ with appropriate vertical arrows we obtain a commutative diagram
$$
\xymatrix{
\cscg\ar[r]^{\SSf\quad}\ar@<-1ex>[d]_{\Lie}&\cscs\ar@<1ex>[l]^{\G}\ar[d]^{\Lie_{\rm lts}}\\
\LAfd\ar[u]\ar[r]^{\R}&\LTSfd\ar@<1ex>[u]\ar@<1ex>[l]^{\A}
}
$$
where $\G$ is an imbedding of $\cscs$ into $\cscg$. A description of the image of $\G$ will follow from

\begin{lm}
Let $(G,\sigma)$ be a connected Lie group with involution and let $g=g_0\oplus g_1$ be the corresponding $\Z_2$-graded Lie algebra. The group $G$ is generated by $G_{\sigma}$ if and only if the algebra $g$ is generated by $g_1$.
\end{lm}

\begin{proof} If $g$ is generated by $g_1$, then it follows from \cite[Lemma~I.3.2]{L} that $G$ is generated by $\{\exp(x):x\in g_1\}$. Besides, $\exp(x)=\exp(\frac{1}{2}x)\sigma(\exp(\frac{1}{2}x))^{-1}\in G_{\sigma}$ for any $x\in g_1$. So, $G$ is generated by $G_{\sigma}$.

To prove the converse we consider the subalgebra $h$ generated by $g_1$. It is easy to see that $h$ is an ideal, so the connected subgroup $H$ corresponding to $h$ is normal. It suffices to prove that $H$ contains $G_{\sigma}$.

There exists a neighborhood $U$ of $0$ in $g$, which is diffeomorphic to a neighborhood $V$ of identity in $G$ via the exponential map $\exp:g\to G$. Also,  there exists a neighborhood $\widetilde U\subseteq U$ of $0$ such that for $W=\exp(\widetilde U)$ one has $WW\subseteq V$, $W^{-1}\subseteq W$, and $W^{\sigma}\subseteq W$.

Since $W$ is a neighborhood of the identity, $G$ is generated by $W=\exp(\widetilde U)$. Thus an arbitrary element $X=Y\sigma(Y)^{-1}\in G_{\sigma}$ can be represented as
\beq
X=\exp(x_1)...\exp(x_n)\exp(-x_n^{\Lie\sigma})...\exp(-x_1^{\Lie \sigma}) \label{XinH}
\eeq
where each $x_i\in \widetilde U$. We use the induction on $n$ to prove that $X\in H$.

For $n=1$, $\exp(x_1)\exp(-x_1^{\Lie \sigma})\in W\sigma(W)^{-1}\subseteq V$ and hence $\exp(x_1)\exp(-x_1^{\Lie \sigma})=\exp(y)$ for some $y\in U$. Moreover, $\exp(y^{\Lie \sigma})=\sigma(\exp(y))=\sigma(\exp(x_1)\exp(-x_1^{\Lie \sigma}))=(\exp(x_1)\exp(-x_1^{\Lie \sigma}))^{-1}
=\exp(-y)$. It follows that $y\in g_1$ and $\exp(x_1)\exp(-x_1^{\Lie \sigma})=\exp(y)\in H$.

Now we consider~(\ref{XinH}) for an arbitrary $n$. We have established that $\exp(x_n)\exp(-x_n^{\Lie \sigma})=\exp(z)\in H$. Consequently,
\begin{equation*}
\begin{split}
X=\exp(x_1)...\exp(x_n)\exp(-x_n^{\Lie \sigma})...\exp(-x_1^{\Lie \sigma})\\
=\exp(x_1)...\exp(x_{n-1})\exp(-x_{n-1}^{\Lie \sigma})...\exp(-x_1^{\Lie \sigma})\\
\times\left[\left(\exp(-x_{n-1}^{\Lie \sigma})...\exp(-x_1^{\Lie \sigma})\right)^{-1}
\exp(z)\left(\exp(-x_{n-1}^{\Lie \sigma})...\exp(-x_1^{\Lie \sigma})\right)\right]\in H.
\end{split}
\end{equation*}
due to the inductive step and the fact that $H$ is normal.
\end{proof}

According to the classical result by Chevaley and Eilenberg the cohomology groups $H^*_{\rm inv}(G)$ of a Lie group $G$ obtained using the left invariant forms can be identified with the cohomology groups of the corresponding Lie algebra $g=\Lie G$. Via this identification we can define the cohomology groups of $(G,\sigma)$ to be $H^*_{\rm gr}(g)$ for the grading of $g$ determined by $\Lie \sigma$. Moreover, since the base ring considered here is a field of characteristic 0 and the vector spaces are finite dimensional, the condition~(ii) of Theorem A is equivalent to $H^*_{\rm gr}(g)=0$. Therefore Theorem A and the lemma above imply Theorem B.

\noindent{\bf Acknowledgement.} I am indebted to John R. Faulkner who brought my attention to the class of Lie triple systems.

\bigskip

\noindent College of Charleston, \\
Charleston, SC 29424, U.S.A.\\
e-mail: smirnov@cofc.edu

\end{document}